\documentclass[11pt,twoside]{article}
\usepackage{amssymb,amsmath}
\usepackage{graphics}
\newcommand{\dl}[1]{{\bf Theorem{#1.}}}
\newcommand{\yl}[1]{{\bf Lemma{#1.}}}
\newcommand{\re}[1]{{\bf Remark{#1.}}}

\newcommand{\de}[1]{{\bf Definition{#1.}}}
\newcommand{\zb}{\hfill$\blacksquare$}
\newcommand{\la}[1]{\label{#1}}
\newcommand{\rf}[1]{(\ref{#1})}

\newcommand{\zm}{{\bf Proof.}}

\begin{document}

\thispagestyle{empty}\setcounter{page}{1}

\begin{center}
{\large An Analysis of the Weak Finite
Element Method for Convection-Diffusion\\ Equations}\\[0.3cm]
\footnotetext{{\em Correspondence to}: Tie Zhang, Department of
Mathematics, Northeastern University,
Shenyang, 110004, China (e-mail: ztmath@ 163.com)\\
Contract grant sponsor: National Natural Science Funds of China:
No. 11371081; and the State Key Laboratory of Synthetical
Automation for Process Industries Fundamental Research Funds, No.
2013ZCX02.}
{Tie Zhang and Yanli Chen\\[0.2cm]}
{\em\small Department
of Mathematics and the State Key Laboratory of Synthetical
Automation for Process Industries, Northeastern University,
Shenyang
110004, China}\\[0.4cm]
\end{center}
We study the weak finite element method solving
convection-diffusion equations. A weak finite element scheme is
presented based on a spacial variational form. We established a
weak embedding inequality that is very useful in the weak finite
element analysis. The optimal order error estimates are derived in
the discrete $H^1$-norm, the $L_2$-norm and the $L_\infty$-norm,
respectively. In particular, the $H^1$-superconvergence of order
$k+2$ is given under certain condition. Finally, numerical
examples are provided to illustrate
our theoretical analysis. \\[0.3cm]
{\em Keywords: weak finite element method; optimal error estimate;
superconvergence; convection-diffusion equation.}

\section{Introduction}
\setcounter{section}{1}\setcounter{equation}{0} Recently, the weak
Galerkin finite element method attracts much attention in the
field of numerical partial differential equations
\cite{Wang,Wang2,Lin,Lin1,Lin2,Lin3,Chen,Li,Har}. This method is
presented originally by Wang and Ye for solving general elliptic problems
in multi-dimensional domain \cite{Wang}. Since then, some modified
weak Galerkin methods have also been studied, for example, see
\cite{Lin4,Gao,Wang0,Yang}. In general, a weak finite element method can be
considered as an extension of the standard finite element method
where classical derivatives are replaced in the variational
equation by the weakly defined derivatives on discontinuous
functions. The main feature of this method is that it allows the
use of totally discontinuous finite element function and the trace
of finite element function on element edge may be independent with
its value in the interior of element. This feature make this
method possess all advantages of the usual discontinuous Galerkin
(DG) finite element method \cite{Arn,Cock,Zhang}, and it has
higher flexibility than the DG method. The readers are referred to
articles \cite{Wang2,Lin,Cock} for more detailed explanation of
this method and its relation with other finite element methods.

In this paper, we study the weak finite element method for
convection-diffusion equations:
\begin{eqnarray}
-\hbox{div}(A\nabla u)+\boldsymbol{b}\cdot\nabla
u+c\,u=f,\;in\;\;\Omega,\la{1.1}
\end{eqnarray}
where $\Omega\subset R^d$, coefficient matrix $A=(a_{ij})_{d\times
d}$ and $\boldsymbol{b}$ is a vector function.

We first introduce the weak gradient and discrete weak gradient
following the way in \cite{Wang}. Then, we consider how to
discretize problem \rf{1.1} by using weak finite elements. In
order to make the weak finite element equation have a positive
property, we present a spacial weak form for problem \rf{1.1}.
This weak form is different from the conventional one and is very
suitable for the weak finite element discretization. We establish
a discrete embedding inequality on the weak finite element space
which provides a useful tool for the weak finite element analysis.
In analogy to the usual finite element research, we derive the
optimal order error estimates in the discrete $H^1$-norm, the
$L_2$-norm and the $L_\infty$-norm, respectively. In particular,
for the pure elliptic problems in divergence form
($\boldsymbol{b}=0$, $c=0$), we obtain an $O(h^{k+2})$-order
superconvergence estimate for the gradient approximation of the
weak finite element solution, when the $(k,k+1)$-order finite
element polynomial pair (interior and edge of element) is used.
Both our theoretical analysis and numerical experiment show that
this weak finite element method is a high accuracy and efficiency
numerical method.

This paper is organized as follows. In Section 2, we establish the
weak finite element method for problem \rf{1.1}. In Section 3,
some approximation functions are given and the stability of the
weak finite element solution is analyzed. Section 4 is devoted to
the optimal error estimate and superconvergence estimate in
various norms. In Section 5, we discuss how to solve the weak
finite element discrete system of equations and then provide some
numerical examples to illustrate our theoretical analysis.

Throughout this paper, for a real $s$, we adopt the notations
$W^{s,p}(D)$ to indicate the usual Sobolev spaces on domain
$D\subset \Omega$ equipped with the norm $\|\cdot\|_{s,p,D}$ and
semi-norm $|\cdot|_{s,p,D}$, and if $p=2$, we set
$W^{s,p}(D)=H^s(D)$, $\|\cdot\|_{s,p,D}=\|\cdot\|_{s,D}$. When
$D=\Omega$, we omit the index $D$. The notations $(\cdot,\cdot)$
and $\|\cdot\|$ denote the inner product and norm in the space
$L_2(\Omega)$, respectively. We will use letter $C$ to represent a
generic positive constant, independent of the mesh size $h$.

\section{Problem and its weak finite element approximation}
\setcounter{section}{2}\setcounter{equation}{0} Consider the
convection-diffusion equations:
\begin{eqnarray}
\left\{\begin{array}{ll}
-\hbox{div}(A\nabla u)+\boldsymbol{b}\cdot\nabla u+c\,u=f,\;in\;\;\Omega,\\
u=g,\;\;on\;\;\partial\Omega,\label{2.1}
\end{array}
 \right.
\end{eqnarray}
where $\Omega\subset R^d\, (d=2,3)$ is a polygonal or polyhedral
domain with boundary $\partial\Omega$, coefficient matrix
$A=(a_{ij})_{d\times d}$ is uniformly positive definite in
$\Omega$, i.e., there exists a positive constant $a_0$ such that
\begin{equation}
a_0\xi^T\xi\leq\xi^TA(x)\xi,\;\forall\,\xi\in R^d,\;x\in
\Omega.\label{2.2}
\end{equation}
We assume that $a_{ij}(x)\in [W^{1,\infty}(\Omega)]^{d\times d}$,
$\boldsymbol{b}\in [W^{1,\infty}(\Omega)]^{d}$ and $c(x)\in
L_\infty(\Omega)$. As usual, we further assume that
\begin{equation}
c(x)-\frac{1}{2}\hbox{div}\boldsymbol{b}(x)\geq 0,\;x\in
\Omega.\la{2.3}
\end{equation}

Let $T_h=\bigcup\lbrace K \rbrace$ be a regular triangulation of
domain $\Omega$ so that $\overline{\Omega}=\bigcup_{K\in
T_h}\{K\}$, where the mesh size $h=max \,h_K$, $h_K$ is the
diameter of element $K$.

In order to define the weak finite element approximation to
problem \rf{2.1}, we first need to introduce the concepts of weak
derivative and discrete weak derivative, which are originally
presented in \cite{Wang}.

Let element $K\in T_h$, denote the interior and boundary of $K$ by
$K^0$ and $\partial K$, respectively. A weak function on element
$K$ refers to a function $v=\{ v^0,v^b\}$ with $v^0=v|_{K^0}\in
L_2(K^0),\,v^b=v|_{\partial K}\in L_2(\partial K)$. Note that for
a weak function $v=\{v^0,v^b\}$, $v^b$ may not be necessarily the
trace of $v^0$
on $\partial K$.\\
\de{ 2.1}\quad The weak derivative $\partial^w_{x_i}v$ of a weak
function $v$ with respect to variable $x_i$ is defined as a linear
functional in the dual space $H^{-1}(K)$ whose action on each
$q\in H^1(K)$ is given by
\begin{equation}
\int_{K}\partial^w_{x_i}vqdx=-\int_{K}v^0\partial_{x_i}qdKdx+\int_{\partial
K}v^bq\cos\theta_ids,\;\forall\, q\in H^1(K),\,1\leq i\leq
d.\la{2.5}
\end{equation}
where $n=(\cos\theta_1,\dots,\cos\theta_d)^T$ is the outward unit
normal vector on $\partial K$.

Obviously, as a bounded linear functional on $H^1(K)$,
$\partial^w_{x_i}v$ is well defined. Moreover, for $v\in H^1(K)$,
if we consider $v$ as a weak function with components
$v^0=v|_{K^0}$, $v^b=v|_{\partial K}$, then by the Green formula,
we have for $q\in H^1(K)$ that
\begin{equation}
\int_{K}\partial_{x_i}vqdx=-\int_{K}v\partial_{x_i}qdx+\int_{\partial
K}vq\cos\theta_ids=-\int_{K}v^0\partial_{x_i}qdx+\int_{\partial
K}v^bq\cos\theta_ids,\la{2.6}
\end{equation}
which implies that $\partial^w_{x_i}v=\partial_{x_i}v$ is the
usual derivative of function $v\in H^1(K)$.

According to Definition 2.1, the weak gradient $\nabla_w v$ of a
weak function $v$ should be such that $\nabla_w
v=(\partial^w_{x_1}v,\dots,\partial^w_{x_d}v)^T\in [H^{-1}(K)]^d$
satisfies
\begin{equation}
\int_{K}\nabla_{w}v\cdot\boldsymbol{q}dx=-\int_Kv^0\hbox{div}\boldsymbol{q}dx+\int_{\partial
K}v^b\boldsymbol{q}\cdot nds,\,\forall\, \boldsymbol{q}\in
[H^{1}(K)]^d.\la{2.7}
\end{equation}

Below we introduce the discrete weak gradient which is actually
used in our analysis.

For any non-negative integer $l\geq 0$, let $P_l(D)$ be the space
composed of all polynomials on set $D$ with degree no more than
$l$. Introduce the discrete weak function space on $K$
\begin{equation}
W(K;k,r)=\{\,v=(v^0,v^b): v^0\in P_k(K^0),\,v^b\in P_{r}(\partial
K)\,\}.\la{2.8}
\end{equation}
\de{ 2.2}\quad For $v\in W(K,k,r)$, the discrete weak derivative
$\partial^w_{x_i,r}v\in P_r(K)$ is defined as the unique solution
of equation:
\begin{equation}
\int_{K}\partial^w_{x_i,r}vqdx=-\int_{K}v^0qdx+\int_{\partial
K}v^bq\cos\theta_ids,\;\forall\, q\in P_r(K),\,1\leq i\leq
d.\la{2.9}
\end{equation}

According to Definition 2.2, for weak function $v\in W(K,k,r)$,
its discrete weak gradient
$\nabla_{w,r}v=(\partial^w_{x_1,r}v,\dots,\partial^w_{x_d,r}v)^T\in
[P_r(K]^d$ is the unique solution of equation:
\begin{equation}
\int_{K}\nabla_{w,r}v\cdot\boldsymbol{q}dx=-\int_Kv^0\hbox{div}\boldsymbol{q}dx+\int_{\partial
K}v^b\boldsymbol{q}\cdot nds,\,\forall\, \boldsymbol{q}\in
[P_r(K)]^d.\la{2.10}
\end{equation}
\re{ 2.1}\quad{\em We here first define the (discrete) weak
derivative, and then the (discrete) weak gradient follows
naturally. This method of defining (discrete) weak gradient is
slightly different from that in \cite{Wang}, in which the
(discrete) weak gradient is defined solely.}

From \rf{2.7} and \rf{2.10}, we have
$$
\int_{K}(\nabla_wv-\nabla_{w,r}v)\cdot\boldsymbol{q}dx=0,\;\forall\,
\boldsymbol{q}\in [P_r(K)]^d.
$$
This shows that $\nabla_{w,r}v$ is a discrete approximation of
$\nabla_wv$ in $[P_r(K)]^d$. In particular, if $v\in H^1(K)$, we
have from \rf{2.6} and \rf{2.10} that
$$
\int_{K}(\nabla
v-\nabla_{w,r}v)\cdot\boldsymbol{q}dx=0,\;\forall\,
\boldsymbol{q}\in [P_r(K)]^d.
$$
That is, $\nabla_{w,r}v$ is the $L_2$ projection of $\nabla v$ in
$[P_r(K)]^d$ if $v\in H^1(K)$.

We have introduced the weak derivative (gradient) and discrete
weak derivative (gradient), but only the discrete weak gradient
given in \rf{2.10} will be used throughout this paper. The others
also should be useful in the study of numerical partial
differential equations.

A important property of $\nabla_{w,r}v$ can be stated as follows,
see \cite[Lemma 5.1]{Wang}.\\
\yl{ 2.1}\quad{\em Let $v=\{ v^0,v^b\}\in W(K,k,r)$ be a weak
function and $r>k$. Then, $\nabla_{w,r}v=0$ on $K$ if and only if
$v=constant$, that is, $v^0=v^b=constant$ on $K$.}

Now, we construct the weak finite element space. Denote by
$\mathcal{E}_h^0=\bigcup\{\,e\in
\partial K\setminus\partial\Omega:\,K\in \mathcal{T}_h\}$ the
union of all boundary faces or edges ($d=2$) of elements in $T_h$
that are not contained in $\partial\Omega$. Let $K_1$ and $K_2$ be
two adjacent elements with the common face $e=\partial
K_1\bigcap\partial K_2$, and $n_1$ and $n_2$ are the outward unit normal
vectors on $\partial K_1$ and $\partial K_2$, respectively. For weak function
defined on $T_h$, set
$v|_{K_1}=\{ v^0_1,v^b_1\},\,v|_{K_2}=\{v^0_2,v^b_2\}$. We define
the jump of weak function $v$ on $e$ by
$$
[v]_{e}=v^b_1n_1+v^b_2n_2=(v_1^b-v_2^b)n_1,\;e\in \mathcal{E}_h^0.
$$
Then, weak function $v$ is single value on $e$ if and only if
$[v]_e=0$. The weak finite element space is now defined by
\begin{eqnarray*}
&&S_h=\{v:\;v|_K\in W(K,k,r),\,K\in
T_h,\,[v]_e=0,\;e\in \mathcal{E}^0_h\},\\
&&S_h^0=\{v:\;v\in S_h,\,v|_{\partial\Omega}=0\}.
\end{eqnarray*}

In order to define the weak finite element approximation to
problem \rf{2.1}, we need to derive a spacial weak form for
problem \rf{2.1}. From the differential formula
\begin{equation}
\bold{b}\cdot\nabla u=\frac{1}{2}\boldsymbol{b}\cdot\nabla
u+\frac{1}{2}\hbox{div}(\boldsymbol{b}\,u)
-\frac{1}{2}\hbox{div}\boldsymbol{b}\,u,\la{2.11}
\end{equation}
and the Green formula, we see that a weak form for problem
\rf{2.1} is to find $u\in H^1(\Omega)$ $u|_{\partial\Omega}=g$
such that
\begin{equation}
(A\nabla u,\nabla v)+\frac{1}{2}(\boldsymbol{b}\cdot\nabla
u,v)-\frac{1}{2}(u,\boldsymbol{b}\cdot\nabla
v)+(c_bu,v)=(f,v),\;\forall\,v\in H^1_0(\Omega),\la{2.12}
\end{equation}
where $c_b=c-\frac{1}{2}\hbox{div}\boldsymbol{b}\geq 0$. Denote
the discrete $L_2$ inner product and norm by
$$
(u,v)_h=\sum_{K\in T_h}(u,v)_{K}=\sum_{K\in
T_h}\int_Ku\,vdx,\;\;\;\;\|u\|_h^2=(u,u)_h.
$$
Motivated by weak form \rf{2.12}, we define the weak finite
element approximation of problem \rf{2.1} by finding $u_h\in
S_h,\,u_h|_{\partial\Omega}=g_h$ such that
\begin{equation}
a_h(u_h,v)=(f,v^0),\,\forall\,v\in S_h^0,\la{2.13}
\end{equation}
where $g_h$ is a proper approximation of function $g$ and the
bilinear form
\begin{equation}
a_h(u,v)=(A\nabla_{w,r} u,\nabla_{w,r}
v)_h+\frac{1}{2}(\boldsymbol{b}\cdot\nabla_{w,r}
u,v^0)_h-\frac{1}{2}(u^0,\boldsymbol{b}\cdot\nabla_{w,r}v)_h+(c_bu^0,v^0).\la{2.14}
\end{equation}
Bilinear form $a_h(u,v)$ is not based on the conventional one:
$$
a(u,v)=(A\nabla u,\nabla v)+(\boldsymbol{b}\cdot\nabla
u,v)+(cu,v).
$$
The advantage of our bilinear form is that it always is positive
definite on the weak function space $S_h\times S_h$, and the
conventional one is not, since the integration by parts does not
hold on weak function space $S_h$ or $S_h^0$.\\
\dl{ 2.1}\quad{\em Let $r>k$. Then, the solution of weak finite
element equation \rf{2.13} uniquely exists. }

\zm\quad Since equation \rf{2.13} is essentially a linear system
of equations, we only need to prove the uniqueness. Let $f=g_h=0$,
we need to prove $u_h=0$. Taking $v=u_h\in S_h^0$ in \rf{2.13}, we
obtain
$$
a_0\|\nabla_{w,r}u_h\|_h^2\leq a_h(u_h,u_h)=0.
$$
This implies that $\nabla_{w,r}u_h=0$ on $T_h$. Thus, from Lemma
2.1, we know that $u_h$ is a piecewise constant on $T_h$. Since
$[u_h]_e=0$ and $u_h|_{\partial\Omega}=0$, so we have $u_h=0$. \zb

\section{Projection and approximation}
\setcounter{section}{3}\setcounter{equation}{0} In this section,
we give some projections and approximation properties which will
be used in next section.

In order to balance the approximation accuracy between spaces
$S_h$ and $P_r(K)$ used to compute $\nabla_{w,r}u_h$, from now on,
we always set the index $r=k+1$ in \rf{2.8}--\rf{2.9}. The other
choice of weak finite element space can be found in
\cite{Wang,Lin}.

For $l\geq 0$, let $P_h^l$ is the local $L_2$ projection operator,
restricted on each element $K$, $P_h^l:\,u\in L_2(K)\rightarrow
P_h^lu\in P_l(K)$ such that
\begin{equation}
(u-P_h^lu,q)_{K}=0,\;\forall\,q\in P_l(K),\,K\in T_h.\la{3.1}
\end{equation}
By the Bramble-Hilbert lemma, it is easy to prove that (see
\cite{Zhang})
\begin{equation}
\|u-P_h^lu\|_{0,K}\leq Ch_K^{s}\|u\|_{s,K},\;0\leq s\leq
l+1.\la{3.2}
\end{equation}
We now define a projection operator $Q_h:\,u\in
H^1(\Omega)\rightarrow Q_hu\in S_h$ such that
\begin{equation}
Q_hu|_{K}=\{Q^0_hu,Q^b_hu\}=\{P_h^ku,P_{\partial
K}^{k+1}u^b\},\,K\in T_h,\la{3.3}
\end{equation}
where $P_{\partial K}^{k+1}$ is the $L_2$ projection operator in
space $P_{k+1}(\partial K)$.\\
\yl{ 3.1}$^{\cite{Wang}}$\quad{\em Let $u\in
H^{1+s}(\Omega),\,s\geq 0$. Then, $Q_hu$ has the following
approximation properties}
\begin{eqnarray}
&&\|u-Q_h^0u\|_{0,K}\leq Ch_K^{s}\|u\|_{s,K},\,0\leq s\leq
k+1,\;K\in T_h,\la{3.4}\\
&&\|\nabla_{w,r}Q_hu-\nabla u\|_{0,K}\leq
Ch^s_K\|u\|_{1+s,K},\,0\leq s\leq k+2,\;K\in T_h.\la{3.5}
\end{eqnarray}

\zm\quad Since $Q_h^0u=P_h^ku$, then estimate \rf{3.4} follows
from \rf{3.2}. Furthermore, from \rf{2.10} and the definition of
$Q_hu$, we have
\begin{eqnarray*}
&&\int_{K}\nabla_{w,r}Q_hu\cdot\boldsymbol{q}dx=-\int_{K}Q_h^0u\hbox{div}\boldsymbol{q}dx+\int_{\partial
K}Q^b_hu \boldsymbol{q}\cdot nds\\
&=&-\int_{K}u\hbox{div}\boldsymbol{q}dx+\int_{\partial
K}u\boldsymbol{q}\cdot nds=\int_{K}\nabla
u\boldsymbol{q}dx,\;\forall\,\boldsymbol{q}\in [P_{r}(K)]^2.
\end{eqnarray*}
This implies $\nabla_{w,r}Q_hu=P_h^r\nabla u$ and estimate
\rf{3.5} holds, noting that $r=k+1$.\zb

For the error analysis, we still need to introduce a special
projection function \cite{Brezzi}. For simplifying, we only
consider the case of two-dimensional domain ($d=2$).

Let $e_i$ and $\lambda_i$ ($1\leq i\leq 3$) are the edge and
barycenter coordinate of $K$, respectively. For function $\phi$,
${\bf curl}\phi=(\partial_{x_2}\phi,-\partial_{x_1}\phi)^T$. Let
Space $P_{k+2}^0(K) =\{\,p\in P_{k+2}(K): p|_{\partial
K}=0\}=\lambda_1\lambda_2\lambda_3P_{k-1}(K)$ and $
H(\hbox{div};\Omega)=\{\boldsymbol{u}\in [L_2(\Omega)]^2:
\hbox{div}\boldsymbol{u}\in L_2(\Omega)\}$.

Define the projection operator $\pi_h:
H(\hbox{div};\Omega)\rightarrow H(\hbox{div};\Omega)$, restricted
on $K\in T_h$, $\pi_h\boldsymbol{u}\in [P_{k+1}(K)]^2$ satisfies
\begin{eqnarray}
&&(\boldsymbol{u}-\pi_h\boldsymbol{u},\nabla
q)_K=0,\;\forall\,q\in
P_k(K),\la{3.6}\\
&&\int_{e_i}(\boldsymbol{u}-\pi_h\boldsymbol{u})\cdot n
qds=0,\;\forall\,q\in P_{k+1}(e_i),\,i=1,2,3,\la{3.7}\\
&&(\boldsymbol{u}-\pi_h\boldsymbol{u},{\bf
curl}\,q)_K=0,\;\forall\, q\in P^0_{k+2}(K).\la{3.8}
\end{eqnarray}
Some properties of projection $\pi_h\boldsymbol{u}$ had been
discussed in \cite{Brezzi}, we here give a more detailed analysis
for our argument requirement.\\
\yl{ 3.2}\quad{\em For $\boldsymbol{u}\in H(\hbox{div}\,;\Omega)$,
the projection $\pi_h\boldsymbol{u}$ uniquely exists and
satisfies}
\begin{eqnarray}
&&(\hbox{div}(\boldsymbol{u}-\pi_h\boldsymbol{u}),q)_K=0,\;\forall\,q\in
P_{k}(K),\;K\in T_h.\la{3.9}
\end{eqnarray}
Furthermore, if $u\in H^{1+s}(\Omega),\,s\geq 0$, then
\begin{eqnarray}
&&\|\pi_h\boldsymbol{u}\|_{0,K}\leq
C\|\boldsymbol{u}\|_{1,K},\;K\in
T_h,\la{3.10}\\
&&\|\boldsymbol{u}-\pi_h\boldsymbol{u}\|_{0,K}\leq
Ch_K^s\|\boldsymbol{u}\|_{s,K},\;1\leq s\leq k+2,\;K\in
T_h.\la{3.11}
\end{eqnarray}

\zm\quad We first prove the unique existence of
$\pi_h\boldsymbol{u}$. Since the number of dimensions (noting that
\rf{3.6} is trivial for $q=constant$):
$$
\hbox{dim}(P_k(K))-1+3\,\hbox{dim}(P_{k+1}(e_i))+\hbox{dim}(P^0_{k+2}(K))=2\,\hbox{dim}(P_{k+1}(K)),
$$
so the linear system of equations \rf{3.6}$\sim$\rf{3.8} is
consistent. Thus, we only need to prove the uniqueness. Assume that
$\boldsymbol{u}=0$ in \rf{3.6}$\sim$\rf{3.8}, we need to prove
$\pi_h\boldsymbol{u}=0$. From \rf{3.6}--\rf{3.7}, we have
$\pi_h\boldsymbol{u}\cdot n=0$ on $\partial K$ and
\begin{eqnarray*}
(\hbox{div}\pi_h\boldsymbol{u},q)_K=-(\pi_h\boldsymbol{u},\nabla
q)_K+\int_{\partial K}\pi_h\boldsymbol{u}\cdot
nqds=0,\;\forall\,q\in P_{k}(K).
\end{eqnarray*}
This implies $\hbox{div}\pi_h\boldsymbol{u}=0$ on $K$. So there
exists a function $\phi\in P_{k+2}(K)$ so that ${\bf
curl}\phi=\pi_h\boldsymbol{u}$ (see \cite{Girault}). Since the
tangential derivative $\partial_\tau\phi=\pi_h\boldsymbol{u}\cdot
n=0$ on $\partial K$, so $\phi=\phi_0=constant$ on $\partial K$.
Let $p=\phi-\phi_0$. Then, $p\in P^0_{k+2}(K)$ and ${\bf
curl}\,p=\pi_h\boldsymbol{u}$. Taking $q=p$ in \rf{3.8}, we obtain
$\|\pi_h\boldsymbol{u}\|_{0,K}=0$ so that $\pi_h\boldsymbol{u}=0$.
Next, we prove \rf{3.9}--\rf{3.11}. Equation \rf{3.9} comes
directly from the Green formula and \rf{3.6}--\rf{3.7}. From the
solution representation of linear system of equations
\rf{3.6}$\sim$\rf{3.8}, it is easy to see that on the reference
element $\hat{K}$,
\begin{equation}
\|\hat{\pi}_h\hat{\boldsymbol{u}}\|_{0,\hat{K}}\leq
\hat{C}(\|\hat{\boldsymbol{u}}\|_{0,\hat{K}}+\|\hat{\boldsymbol{u}}\|_{0,\partial
\hat{K}})\leq\hat{C}(\|\hat{\boldsymbol{u}}\|_{0,\hat{K}}+\|\nabla\hat{\boldsymbol{u}}\|_{0,\hat{K}}),\la{3.12}
\end{equation}
where we have used the trace inequality. Then,
\rf{3.10} follows from \rf{3.12} and a scale argument between $\hat{K}$ and $K$. From
\rf{3.12}, we also obtain
$$
\|\hat{\pi}_h\hat{\boldsymbol{u}}\|_{0,\hat{K}}\leq
\hat{C}\|\hat{\boldsymbol{u}}\|_{s,\hat{K}}, 1\leq s\leq k+2.
$$
Hence, estimate \rf{3.11} can be derived by using the
Bramble-Hilbert lemma.\zb

The following discrete embedding inequality is an analogy of the
Poincar\'e inequality in $H^1_0(\Omega)$.\\
\yl{ 3.3}\quad{\em Let $\Omega$ be a polygonal or polyhedral
domain. Then, for $v\in S_h^0$, there is a positive constant $C_0$
independent of $h$ such that}
\begin{equation}
\|v^0\|\leq C_0\|\nabla_{w,r}v\|_h,\;\forall\, v\in
S_h^0.\la{3.13}
\end{equation}

\zm\quad For $v\in S^0_h$, we first make a smooth domain
$\Omega'\supset \Omega$ ( if $\Omega$ is convex, we may set
$\Omega'=\Omega$) and extend $v^0$ to domain $\Omega'$ by setting
$v^0|_{\Omega'\backslash\Omega}=0$. Then, there exists a function
$w\in H^1_0(\Omega')\bigcap H^2(\Omega')$ such that
\begin{eqnarray*}
-\triangle w=v^0,\;in\,\,\Omega',\;\;\|w\|_{2,\Omega'}\leq C\|v^0\|.
\end{eqnarray*}
Now we set $\boldsymbol{w}=-\nabla w$, then $\boldsymbol{w}\in [H^1(\Omega)]^d$ satisfies
$$
\hbox{div}\boldsymbol{w}=v^0,\;in\,\,\Omega,\;\;\|\boldsymbol{w}\|_1\leq
\|w\|_{2,\Omega'}\leq C\|v^0\|.
$$
Hence, we have from \rf{3.9}, \rf{3.10} and \rf{2.10} that
\begin{eqnarray*}
\|v^0\|^2&=&(\hbox{div}\boldsymbol{w},v^0)=(\hbox{div}\pi_h\boldsymbol{w},v^0)\\
&=&\sum_{K\in T_h}\big(-\int_K\nabla_{w,r}v\cdot
\pi_h\boldsymbol{w}dx+\int_{\partial
K}v^b\pi_h\boldsymbol{w}\cdot nds\big)\\
&=&\sum_{K\in
T_h}-\int_K\nabla_{w,r}v\cdot\pi_h\boldsymbol{w}dx\leq
\|\nabla_{w,r}v\|_h\|\pi_h\boldsymbol{w}\|_h\\
&\leq& C\|\nabla_{w,r}v\|_h\|\boldsymbol{w}\|_1\leq
C\|\nabla_{w,r}v\|_h\|v^0\|,
\end{eqnarray*}
where we have used the fact that $[v]_e=0$ and
\begin{equation}
\sum_{K\in T_h}\int_{\partial K}v^b\pi_h\boldsymbol{w}\cdot
nds=\sum_{K\in T_h}\int_{\partial K}v^b\boldsymbol{w}\cdot
nds=\sum_{e\in \mathcal{E}^0_h}\int_{e}[v]_e\boldsymbol{w}\cdot
nds=0.\la{3.14}
\end{equation}
The proof is completed.\zb

A direct application of Lemma 3.3 is the stability estimate of
weak finite element solution $u_h$.\\
\yl{ 3.4}\quad{\em Let $u_h\in S_h$ be the solution of problem
\rf{2.13}, $g_h=Q_h^bg$ and $u\in H^1(\Omega)$ the solution of problem
\rf{2.1} with $f\in L_2(\Omega)$ and $g\in
H^{\frac{1}{2}}(\partial\Omega)$. Then we have}
$$
\|u^0_h\|+\|\nabla_{w,r}u_h\|_h\leq
C(\|f\|+\|g\|_{\frac{1}{2},\partial\Omega}).
$$

\zm\quad Let $e_h=u_h-Q_hu$. Then, from \rf{2.13}, we see that
$e_h\in S_h^0$ satisfies
$$
a_h(e_h,v)=(f,v^0)-a_h(Q_hu,v),\,v\in S_h^0.
$$
Taking $v=e_h$ and noting that $Q_h^0u=P_h^ku$ and
$\nabla_{w,r}Q_hu=P_h^r\nabla u$, we have
\begin{eqnarray*}
\|\nabla_{w,r}e_h\|^2_h\leq
\|f\|\|e_h^0\|+(|A|_\infty+|\boldsymbol{b}|_\infty+|c_b|_\infty)(\|\nabla
u\|+\|u\|)(\|\nabla_{w,r}e_h\|_h+\|e^0_h\|).
\end{eqnarray*}
Using embedding inequality \rf{3.13} and the a priori estimate for
elliptic problem \rf{2.1}, the proof is completed.\zb
\section{Error analysis}
\setcounter{section}{4}\setcounter{equation}{0} In this section,
we do the error analysis for the weak finite element method
\rf{2.13}. We will see that the weak finite element method
possesses the same (or better) convergence order as that of the
conventional finite element method.

In following error analysis, we always assume that the data
$A, \boldsymbol{b}$ and $c$ in problem \rf{2.1} is smooth enough
for our argument.\\
\yl{ 4.1}\quad{\em Let $u\in H^2(\Omega)$ be the solution of
problem \rf{2.1}. Then we have}
$$
(\pi_h(A\nabla
u),\nabla_{w,r}v)_h+\frac{1}{2}(\boldsymbol{b}\cdot\nabla
u,v^0)_h-\frac{1}{2}(\pi_h(\boldsymbol{b}u),\nabla_{w,r}v)_h+(c_bu,v^0)=(f,v^0),\;v\in
S_h^0.
$$

\zm\quad Let $\boldsymbol{w}\in [H^1(\Omega)]^d$. From Lemma 3.2
and \rf{3.14}, we have
\begin{equation}
-(\hbox{div}\boldsymbol{w},v^0)_h=-(\hbox{div}\pi_h\boldsymbol{w},v^0)_h
=(\pi_h\boldsymbol{w},\nabla_{w,r}v)_h,\;v\in S_h^0.\la{4.1}
\end{equation}
Next, from equations \rf{2.1} and \rf{2.11}, we obtain
$$
-(\hbox{div}(A\nabla
u),v^0)_h+\frac{1}{2}(\boldsymbol{b}\cdot\nabla
u,v^0)_h-\frac{1}{2}(\hbox{div}(\boldsymbol{b}u),v^0)_h+(c_bu,v^0)=(f,v^0),\;v\in
S_h^0.
$$
Together with \rf{4.1} in which setting $\boldsymbol{w}=A\nabla u$
and $\boldsymbol{w}=\boldsymbol{b}u$, respectively, we arrive at the conclusion.
\zb

We first give an abstract error estimate for $u_h-Q_hu$ in the discrete $H^1$-norm.\\
\dl{ 4.1}\quad{\em Let $u$ and $u_h$ be the solutions of problems
\rf{2.1} and \rf{2.13}, respectively, $u\in H^2(\Omega)$ and
$g_h=Q_h^bg$. Then, we have}
\begin{eqnarray}
a_0\|\nabla_{w,r}u_h-\nabla_{w,r}Q_hu\|_h&\leq&
\|A\nabla_{w,r}Q_hu-\pi_h(A\nabla
u)\|_h+C_0\|\boldsymbol{b}\cdot(\nabla_{w,r}Q_hu-\nabla u)\|_h\nonumber\\
&&+\|\boldsymbol{b}Q_h^0u-\pi_h(\boldsymbol{b}u)\|_h+C_0\|c_b(Q_h^0u-u)\|_h.\la{4.2}
\end{eqnarray}

\zm\quad From Lemma 4.1, we see that $Q_hu$ satisfies the equation
\begin{eqnarray}
a_h(Q_hu,v)&=&(f,v^0)+(A\nabla_{w,r}Q_hu-\pi_h(A\nabla
u),\nabla_{w,r}v)_h\nonumber\\
&&+\frac{1}{2}(\boldsymbol{b}\cdot(\nabla_{w,r}Q_hu-\nabla
u),v^0)_h
-\frac{1}{2}(\boldsymbol{b}Q_h^0u-\pi_h(\boldsymbol{b}u),\nabla_{w,r}v)_h\nonumber\\
&&+(c_b(Q_h^0u-u),v^0)_h,\;\forall\,v\in S_h^0.\la{4.3}
\end{eqnarray}
Combining this with equation \rf{2.13}, we obtain the error
equation
\begin{eqnarray}
a_h(Q_hu-u_h,v)&=&(A\nabla_{w,r}Q_hu-\pi_h(A\nabla
u),\nabla_{w,r}v)_h\nonumber\\
&&+\frac{1}{2}(\boldsymbol{b}\cdot(\nabla_{w,r}Q_hu-\nabla
u),v^0)_h
-\frac{1}{2}(\boldsymbol{b}Q_h^0u-\pi_h(\boldsymbol{b}u),\nabla_{w,r}v)_h\nonumber\\
&&+(c_b(Q_h^0u-u),v^0)_h,\;\forall\,v\in S_h^0.\la{4.4}
\end{eqnarray}
Taking $v=Q_hu-u_h$ in \rf{4.4} and using embedding inequality \rf{3.13}, we arrive at the conclusion of Theorem 4.1.\zb

By means of Theorem 4.1, we can derive the following error
estimates.\\
\dl{ 4.2}\quad{\em Let $u$ and $u_h$ be the solutions of problems
\rf{2.1} and \rf{2.13}, respectively, $u\in H^{2+s}(\Omega), s\geq
0$, and $g_h=Q_h^bg$. Then we have the optimal order error
estimates
\begin{eqnarray}
\|u_h^0-u\|+\|\nabla_{w,r}u_h-\nabla u\|_h \leq
Ch^{1+s}\|u\|_{2+s},\;0\leq s\leq k.\la{4.5}
\end{eqnarray}
In particular, for the pure elliptic problem in divergence form
($\boldsymbol{b}=0, c=0$) and $u \in H^{3+s}(\Omega),\,s\geq 0$,
we have the superconvergence estimate}
\begin{eqnarray}
\|\nabla_{w,r}u_h-\nabla u\|_h \leq Ch^{2+s}\|u\|_{3+s},\;0\leq
s\leq k.\la{4.6}
\end{eqnarray}

\zm\quad Using the approximation properties \rf{3.4}--\rf{3.5} and
\rf{3.11}, we obtain
\begin{eqnarray*}
&&\|A\nabla_{w,r}Q_hu-\pi_h(A\nabla u)\|_h\\
&\leq& |A|_\infty\|\nabla_{w,r}Q_hu-\nabla u\|_h+\|A\nabla
u-\pi_h(A\nabla
u)\|_h\leq Ch^{1+s}\|u\|_{2+s},\,0\leq s\leq k+1,\\
&&\|\boldsymbol{b}\cdot(\nabla_{w,r}Q_hu-\nabla u)\|_h\leq Ch^{1+s}\|u\|_{2+s},\;0\leq s\leq k+1,\\
&&\|\boldsymbol{b}Q_h^0u-\pi_h(\boldsymbol{b}u)\|_h+\|c_b(Q_h^0u-u)\|_h\leq
Ch^{1+s}\|u\|_{1+s},\,0\leq s\leq k.
\end{eqnarray*}
Substituting these estimates into \rf{4.2}, we obtain
\begin{equation}
\|\nabla_{w,r}u_h-\nabla_{w,r}Q_hu\|_h\leq
Ch^{1+s}\|u\|_{2+s},\,0\leq s\leq k.\la{4.7}
\end{equation}
Hence, using the triangle inequality
$$
\|\nabla_{w,r}u_h-\nabla u\|_h\leq
\|\nabla_{w,r}u_h-\nabla_{w,r}Q_hu\|+\|\nabla_{w,r}Q_hu-\nabla
u\|_h
$$
and approximation property \rf{3.5}, estimate \rf{4.5} is derived
for the discrete $H^1$-norm. Since
$$
\|u_h^0-u\|\leq \|u^0_h-Q_h^0u\|+\|Q_h^0u-u\|,
$$
then the $L_2$-error estimate follows from the discrete embedding
inequality \rf{3.13}, estimates \rf{4.7} and \rf{3.4}.

Furthermore, if $\boldsymbol{b}=0, c=0$, from Theorem 4.1, we have
$$
a_0\|\nabla_{w,r}u_h-\nabla_{w,r}Q_hu\|_h\leq
\|A\nabla_{w,r}Q_hu-\pi_h(A\nabla u)\|_h.
$$
Then, the superconvergence estimate \rf{4.6} can be derived by
using approximation properties \rf{3.5} and \rf{3.11} and the
triangle inequality.\zb

Theorem 4.2 shows that the weak finite element method is a high
accuracy numerical method, in particular, in the gradient
approximation.

From Theory 4.2, we see that, for $k$-order finite element, weak
finite element method usually has higher accuracy than other
finite element methods in the gradient approximation. The reason
is that the discrete weak gradient is computed by using
$(k+1)$-order polynomial. This action will add some computation
expense, but all additional computations are implemented locally,
in the element level, see \rf{2.10} and Section 5.

Below we give a superclose estimate for error $Q^0_hu-u_h^0$. To
this end, we assume problem \rf{2.1} has the $H^2$-regularity and consider the auxiliary problem: $w\in
H^1_0(\Omega)\bigcap H^2(\Omega)$ satisfies
\begin{eqnarray}
-\hbox{div}(A\nabla w)-\boldsymbol{b}\cdot\nabla
w+c\,w=Q_h^0u-u_h,\;in\;\;\Omega,\;\|w\|_2\leq
C\|Q_h^0u-u_h^0\|.\label{4.8}
\end{eqnarray}
From the argument of Lemma 4.2, we know that $w$ satisfies
equation:
\begin{eqnarray}
&&(\pi_h(A\nabla w),\nabla_{w,r}v)_h-
\frac{1}{2}(\boldsymbol{b}\cdot\nabla
w,v^0)_h+\frac{1}{2}(\pi_h(\boldsymbol{b}w),\nabla_{w,r}v)_h+(c_bw,v^0)\nonumber\\
&=&(Q_h^0u-u_h^0,v^0),\;\forall\,v\in S_h^0.\la{4.9}
\end{eqnarray}
\dl{ 4.3}\quad{\em Let $u$ and $u_h$ be the solutions of problems
\rf{2.1} and \rf{2.13}, respectively, $u\in H^{2+s}(\Omega), s\geq
0, g_h=Q_h^bg$. Then we have the following superclose estimate}
\begin{eqnarray}
\|Q_h^0u-u_h^0\|\leq Ch^{2+s}\|u\|_{2+s},\;0\leq s\leq k,\la{4.10}
\end{eqnarray}
and the optimal $L_2$-error estimate
\begin{eqnarray}
\|u_h^0-u\|\leq Ch^{k+1}\|u\|_{k+1},\;k\geq 1\,.\la{4.10a}
\end{eqnarray}

\zm\quad Taking $v=Q_hu-u_h$ in \rf{4.9}, we have
\begin{eqnarray}
&&\|Q_h^0u-u_h^0\|^2=(\nabla_{w,r}(Q_hu-u_h),\pi_h(A\nabla
w))_h\nonumber\\
&&-\frac{1}{2}(\boldsymbol{b}\cdot\nabla
w,Q_h^0u-u_h^0)_h+\frac{1}{2}(\pi_h(\boldsymbol{b}w),\nabla_{w,r}(Q_hu-u_h))_h+(c_b(Q_h^0u-u_h^0),w)\nonumber\\
&=&(\nabla_{w,r}(Q_hu-u_h),\pi_h(A\nabla
w)-A\nabla_{w,r}Q_hw)_h-\frac{1}{2}(\boldsymbol{b}\cdot(\nabla
w-\nabla_{w,r}Q_hw),Q_h^0u-u_h^0)_h\nonumber\\
&&+\frac{1}{2}(\pi_h(\boldsymbol{b}w)-\boldsymbol{b}Q_h^0w,\nabla_{w,r}(Q_hu-u_h))_h+(c_b(w-Q_h^0w),Q_h^0u-u_h^0)\nonumber\\
&&+a_h(Q_hu-u_h,Q_hw)\nonumber\\
&\leq& Ch\|w\|_2\big(\|\nabla_{w,r}(Q_hu-u_h)\|_h+\|Q_h^0u-u_h^0\|\big)+a_h(Q_hu-u_h,Q_hw)\nonumber\\
&\leq& Ch^{2+s}\|u\|_{2+s}\|w\|_2+a_h(Q_hu-u_h,Q_hw),\,0\leq s\leq k,\la{4.11}
\end{eqnarray}
where we have used embedding inequality \rf{3.13},
estimate \rf{4.7} and the approximation properties of
$\pi_hw,\,\nabla_{w,r}Q_hw$ and $Q^0_hw$. Below we only need to
estimate $a_h(Q_hu-u_h,Q_hw)$.  Using error equation \rf{4.4}, we
have
\begin{eqnarray}
a_h(Q_hu-u_h,Q_hw)&=&(A\nabla_{w,r}Q_hu-\pi_h(A\nabla
u),\nabla_{w,r}Q_hw)_h\nonumber\\
&&+\frac{1}{2}(\boldsymbol{b}\cdot(\nabla_{w,r}Q_hu-\nabla
u),Q^0_hw)_h
-\frac{1}{2}(\boldsymbol{b}Q_h^0u-\pi_h(\boldsymbol{b}u),\nabla_{w,r}Q_hw)_h\nonumber\\
&&+(c_b(Q_h^0u-u),Q^0_hw)_h,\;\forall\,v\in
S_h^0\nonumber\\
&=&E_1+E_2+E_3+E_4.\la{4.12}
\end{eqnarray}
Since $Q_h^0=P_h^k$ and $\nabla_{w,r}Q_h=P_h^r$, then by using
Green's formula and Lemma 3.2, we have
\begin{eqnarray*}
E_1&=&(A\nabla_{w,r}Q_hu-\pi_h(A\nabla u),\nabla_{w,r}Q_hw-\nabla
w)_h+(A\nabla_{w,r}Q_hu-\pi_h(A\nabla u),\nabla w)_h\\
&\leq& Ch^{2+s}\|u\|_{2+s}\|w\|_2+(A\nabla_{w,r}Q_hu-A\nabla
u,\nabla w)_h+(A\nabla u-\pi_h(A\nabla
u),\nabla w)_h\\
&=& Ch^{2+s}\|u\|_{2+s}\|w\|_2+(\nabla_{w,r}Q_hu-\nabla
u,A^T\nabla w-P_h^k(A^T\nabla w))_h\\
&&-(\hbox{div}(A\nabla u-\pi_h(A\nabla u)),w-P_h^kw)_h\leq
Ch^{2+s}\|u\|_{2+s}\|w\|_2,\,0\leq s\leq k.
\end{eqnarray*}
Similarly, we obtain
\begin{eqnarray*}
E_2&=&\frac{1}{2}(\boldsymbol{b}\cdot(\nabla_{w,r}Q_hu-\nabla
u),Q^0_hw)_h=\frac{1}{2}(\boldsymbol{b}\cdot(\nabla_{w,r}Q_hu-\nabla
u),Q^0_hw-w)_h\\
&&+\frac{1}{2}(\nabla_{w,r}Q_hu-\nabla
u,\boldsymbol{b}w-P_h^0(\boldsymbol{b}w))_h\leq
Ch^{2+s}\|u\|_{2+s}\|w\|_1,\,0\leq s\leq k,\\
E_3&=&-\frac{1}{2}(\boldsymbol{b}Q_h^0u-\pi_h(\boldsymbol{b}u),\nabla_{w,r}Q_hw-\nabla
w)_h-\frac{1}{2}(\boldsymbol{b}Q_h^0u-\pi_h(\boldsymbol{b}u),\nabla
w)_h\\
&\leq&
Ch^{2+s}\|u\|_{2+s}\|w\|_2-\frac{1}{2}(\boldsymbol{b}Q_h^0u-\boldsymbol{b}u,\nabla
w)_h-\frac{1}{2}(\boldsymbol{b}u-\pi_h(\boldsymbol{b}u),\nabla w)_h\\
&=&Ch^{2+s}\|u\|_{2+s}\|w\|_2-\frac{1}{2}(Q_h^0u-u,\boldsymbol{b}\nabla
w-P_h^k(\boldsymbol{b}\nabla w))_h\\
&\leq& Ch^{2+s}\|u\|_{2+s}\|w\|_2,\,0\leq s\leq k,\\
E_4&=&(c_b(Q_h^0u-u),Q^0_hw)_h=(c_b(Q_h^0u-u),Q^0_hw-w)_h\\
&&+((Q_h^0u-u),c_bw-P_h^k(c_bw))_h\leq
Ch^{2+s}\|u\|_{1+s}\|w\|_1,\,0\leq s\leq k.
\end{eqnarray*}
Substituting estimates $E_1\sim E_4$ into \rf{4.12} and combining
\rf{4.11}, we arrive at estimate \rf{4.10}, noting that
$\|w\|_2\leq C\|Q_h^0u-u_h^0\|$. Estimate \rf{4.10a} follows from
\rf{4.10} and the triangle inequality. \zb

The difference between estimates \rf{4.5} and \rf{4.10a} is that
for getting the $O(h^{k+1})$-order error estimate in the
$L_2$-norm, the regularity requirement in \rf{4.10a} is optimal
and lower than that in \rf{4.5}.

In order to derive the $L_\infty$-error estimate, we need to
impose the quasi-uniform condition on partition $T_h$ so that the
finite element inverse inequality holds in $S_h$. \\
\dl{ 4.4}\quad{\em Assume that partition $T_h$ is quasi-uniform,
and $u$ and $u_h$ are the solution of problems \rf{2.1} and
\rf{2.13}, respectively, $u\in W^{1+s,\infty}(\Omega)\bigcap
H^{2+s}(\Omega), s\geq 0,\,g_h=Q_h^bg$. Then, we have}
\begin{equation}
\|u-u_h^0\|_{0,\infty}\leq
Ch^{2+s-\frac{d}{2}}(\|u\|_{1+s,\infty}+\|u\|_{2+s}), \,0\leq
s\leq k.\la{4.14}
\end{equation}

\zm\quad From Theorem 4.3 and the finite element inverse
inequality, we have that
$$
\|Q_h^0u-u_h^0\|_{0,\infty}\leq
Ch^{-\frac{d}{2}}\|Q_h^0u-u_h^0\|\leq
Ch^{2+s-\frac{d}{2}}\|u\|_{2+s}.
$$
Hence, by using the approximation property of $Q_h^0u=P_h^ku$, we
obtain
\begin{eqnarray*}
\|u-u_h^0\|_{0,\infty}&\leq&
\|u-Q_h^0u\|_{0,\infty}+\|Q_h^0u-u_h^0\|_{0,\infty}\\
&\leq& Ch^{1+s}\|u\|_{1+s,\infty}+Ch^{2+s-\frac{d}{2}}\|u\|_{2+s}.
\end{eqnarray*}
The proof is completed. \zb

For two-dimensional problem $(d=2)$, Theorem 4.4 gives the optimal
order error estimate in the $L_\infty$-norm.
\section{Numerical experiment}
\setcounter{section}{5} \setcounter{equation}{0} In this section,
we discuss how to solve the weak finite element equation \rf{2.13}
and give some numerical examples to illustrate our theoretical
analysis.

\subsection{Weak finite element linear system of equations}

In order to form the discrete linear system of equations from weak
finite element equation \rf{2.13}, we first introduce the basis
functions of space $S_h^0$. Let $K$ is a element and $e$ is an
edge of $K$. Further let $\{\varphi_{j,K}(x),\,j=1,\dots,N_k\}$ be
the basis functions of space $P_k(K)$,
$\{\varphi_{j,e}(x),\,j=1,\dots,N_e\}$ the basis functions of
space $P_{k+1}(e)$ and $\{\varphi_{j,b}(x),\,j=1,\dots,M_b\}$ be
the basis functions of polynomial set $\{p\in P_{k+1}(\partial
K),\,K\in T_h: p|_e=\varphi_{j,e}(x),\,e\in \mathcal{E}_h^0\}$.
Set the weak basis functions $\psi_{j,K}=\{\varphi_{j,K},0\}$ and
$\psi_{j,b}=\{0,\varphi_{j,b}\}$. Then, we have
$S_h^0=span\{\psi_{1,K},\dots,\psi_{N_k,K},\psi_{1,b},\dots,\psi_{M_b,b},
K\in T_h\}$. By definition \rf{2.10} of discrete weak gradient, we
see that the support set of $\nabla_{w,r}\psi_{j,K}(x)$ is in $K$
and the support set of $\nabla_{w,r}\psi_{j,b}(x)$ is in $K^b=\{K:
\partial K\bigcap supp(\psi_{j,b})\neq\emptyset\}$. Thus, weak finite
element equation \rf{2.13} is equivalent to the following linear
system of equations: $u_h\in S_h$, $u_h|_{\partial\Omega}=g_h$,
such that
\begin{eqnarray}
&&a_K(u_h,v)=(f,v^0)_K,\,v\in\{\psi_{j,K}\},\, K\in T_h,\, \la{5.1}\\
&&a_{K^b}(u_h,v)=0,\,v\in\{\psi_{j,b}\},\,j=1,\dots,M_b,\la{5.2}
\end{eqnarray}
where $a_D(u,v)$ is the restriction of $a_h(u,v)$ on set $D$,
i.e., all integrals in $a_h(u,v)$ are restricted on $D$. Equations
\rf{5.1}--\rf{5.2} form a linear system composed of $N\times
N_k+M_b$ equations with $N\times N_k+M_b$ unknowns, where $N$ is
the total number of elements in $T_h$.

To solve equations \rf{5.1}--\rf{5.2}, we need to design a solver
to compute the discrete weak gradient $\nabla_{w,r}v$ or
$\nabla_{w,r}u_h$. According to \rf{2.10}, for given
$v=\{v^0,v^b\}\in W(K,k,r)$ ($r=k+1$), $\nabla_{w,r}v\in
[P_{k+1}(K)]^d$ can be computed by the following formula
\begin{equation}
M_KV_{w,r}=A_KV^0+\sum_{e\in\partial K}B_{e}V^e, \la{5.3}
\end{equation}
where $V_{w,r}$, $V^0$ and $V^e$ are the vectors associated with
functions $\nabla_{w,r}v\in [P_{k+1}(K)]^d$, $v^0\in P_k(K)$ and
$v^e=v^b|_{e}\in P_{k+1}(e)$, respectively. The matrixes in
\rf{5.3} are as follows
\begin{eqnarray*}
&&M_K=(m_{ij})_{dN_r\times dN_r}, \,\,A_K=(a_{ij})_{dN_r\times
N_k},\,\,B_{e}=(b_{ij})_{dN_r\times N_e},\\
&&m_{ij}=(\overrightarrow{\varphi}_{j,K},\overrightarrow{\varphi}_{i,K})_K,\,\,a_{ij}
=-(\varphi_{j,K},\hbox{div}\overrightarrow{\varphi}_{i,K})_{K},\,
\,b_{ij}=(\varphi_{j,b},\overrightarrow{\varphi}_{i,K}\cdot
n)_{e}.
\end{eqnarray*}
Now, linear system of equations \rf{5.1}--\rf{5.2} can be solved
in the following way. We first use formula \rf{5.3} to derive the
linear relation $\nabla_{w,r}u_h(K)=L(u_h^0(K),u^b_h(\partial
K))$. Then, by substituting $\nabla_{w,r}u_h(K)$ into equations
\rf{5.1}--\rf{5.2}, we can obtain a linear system of equations
that only concerns unknowns $\{u_h^0(K),u_h^b(\partial K)\}, K\in
T_h$. This linear system can be solved by using a proper linear
solver, in which $\nabla_{w,r}v(K)$ is computed by formula
\rf{5.3}.

\subsection{Numerical example}

Let us consider problem \rf{2.1} with the following data:
\begin{equation}
A=(1+x_1x_2)I,\,
\boldsymbol{b}=(1,2)^T,\,c=\sin(x_1x_2),\,u=\sin(x_1\pi)\sin(x_2\pi),\la{5.4}
\end{equation}
and the source term $f =-\hbox{div}(A\nabla
u)+\boldsymbol{b}\cdot\nabla u+cu$, domain $\Omega=(0, 1)^2$.

In the numerical experiments, we first partition $\Omega$ into a
regular triangle mesh $T_h$ with mesh size $h=1/N$. Then, the
refined mesh $T_{h/2}$ is obtained by connecting the midpoint of
each edge of elements in $T_h$ by straight line. Thus, we obtain a
mesh series $T_{h/2^j}, j=0,1,\dots$. We use the polynomial pair
$(P_0(K),P_1(\partial K))$ for space $W(K,k,r)$ and set
$P_r(K)=P_1(K)$ in the weak finite element discretization
\rf{2.13}. We examine the computation error in the discrete
$H^1$-norm, the $L_2$-norm and the $L_\infty$-norm, respectively.
The numerical convergence rate is computed by using the formula
$r=\ln(e_h/e_{\frac{h}{2}})/\ln 2$, where $e_h$ is the computation
error. Table I give the numerical results. We see that the
convergence rates are consistent with or better than the
theoretical prediction. Then, we further examine the
superconvergence of the weak finite element solution (see Theorem
4.2). We take the data as in \rf{5.4} with $\boldsymbol{b}=0$ and
$c=0$. The desired $H^1$-superconvergence is observed from the
numerical results, see Table II. In particular, from Table I-II,
we see that the numerical convergence rates in $L_2$-norm and
$L_\infty$-norm are also superconvergent, although we have no such
conclusion in theory.

\begin{center}
{\small{TABLE I}\quad History of convergence}\\[0.5\baselineskip]
\renewcommand\arraystretch{1.4}
\arrayrulewidth 0.5pt \small
\begin{tabular}{cccc}
\hline
&$\|\nabla_{w,r}u_h-\nabla u\|_h$&\quad$\|u-u_h^0\|$&\quad$\|u_h^0-u\|_\infty$\\
mesh $h$&error \quad\quad rate&\qquad error \qquad\quad rate&\quad error \qquad\quad rate\\
\hline
1/4&0.9329\ \  \quad\ \  - &\quad \ \ 0.870e-1 \qquad\ \ \ -  &1.990e-1\qquad \ \ \ \ \ \ -\\
1/8&4.660e-1\quad  1.0013 &\qquad2.427e-2\qquad 1.8417  &5.525e-2\qquad 1.8489\\
1/16&2.325e-1\quad 1.0033 &\qquad6.246e-3\qquad 1.9579  &1.413e-2\qquad 1.9677\\
1/32&1.162e-1\quad 1.0011 &\qquad1.573e-3\qquad 1.9893  &3.555e-3\qquad 1.9905\\
1/64&5.806e-2\quad 1.0003 &\qquad3.940e-4\qquad 1.9973  &8.906e-4\qquad 1.9969\\
1/128&2.903e-2\quad1.0001 &\qquad0.986e-4\qquad 1.9993  &2.228e-4\qquad 1.9993\\
 \hline
\end{tabular}
\end{center}
\begin{center}
{\small{TABLE II}\quad History of convergence with
$\boldsymbol{b}=0$ and
$c=0$}\\[0.5\baselineskip]
\renewcommand\arraystretch{1.4}
\arrayrulewidth 0.5pt \small
\begin{tabular}{cccc}
\hline
&$\|\nabla_{w,r}u_h-\nabla u\|_h$&\quad$\|u-u_h^0\|$&\quad$\|u_h^0-u\|_\infty$\\
mesh $h$&error \quad\quad rate&\qquad error \qquad\quad rate&\quad error \qquad\quad rate\\
\hline
1/4&1.875e-1\ \ \ \quad\ \  - &\quad \ \ 3.129e-2 \qquad\ \ \ -  &0.928e-1\qquad \ \ \ \ \ \ -\\
1/8&4.896e-2\quad  1.9370 &\qquad8.538e-3\qquad 1.8735  &2.553e-2\qquad 1.8625\\
1/16&1.239e-2\quad 1.9821 &\qquad2.184e-3\qquad 1.9673  &6.533e-3\qquad 1.9664\\
1/32&3.109e-3\quad 1.9948 &\qquad5.490e-4\qquad 1.9917  &1.642e-3\qquad 1.9919\\
1/64&7.782e-4\quad 1.9984 &\qquad1.375e-4\qquad 1.9979  &4.112e-4\qquad 1.9979\\
1/128&1.946e-4\quad1.9995 &\qquad3.437e-5\qquad 1.9995  &1.028e-4\qquad 1.9993\\
 \hline
\end{tabular}
\end{center}

\end{document}